\documentclass[11pt,letterpaper,reqno]{amsart}
\usepackage{amssymb,amscd,amsbsy,mathrsfs}
\newcommand{\s}{\sigma}

\newcommand{\f}{\varphi}

\newcommand{\B}{{\mathscr B}}

\newcommand{\E}{{\mathscr E}}

\newcommand{\F}{{\mathscr F}}

\newcommand{\X}{{\mathscr X}}
\newcommand{\Y}{{\mathscr Y}}

\newcommand{\R}{{\Bbb R}}

\newcommand{\bS}{{\boldsymbol S}}

\newcommand{\rf}[1]{(\ref{#1})}

\newcommand{\df}{\stackrel{\mathrm{def}}{=}}

\newcommand{\supp}{\operatorname{supp}}

\newcommand{\const}{\operatorname{const}}

\newcommand{\eeq}{\end{equation}}
\newcommand{\beq}{\begin{equation}}
\newcommand{\bay}{\begin{eqnarray}}
\newcommand{\ba}{\begin{align*}}
\newcommand{\ea}{\end{align*}}
\newcommand{\ey}{\end{eqnarray}}
\newcommand{\bey}{\begin{eqnarray*}}
\newcommand{\eey}{\end{eqnarray*}}

\newcommand{\be}{\infty}

\newtheorem{thm}{\hspace{\parindent}Theorem}[section]

\newtheorem{cor}[thm]{\hspace{\parindent}Corollary}

\usepackage{amsmath}
\usepackage{amsfonts}
\usepackage{amssymb}
\usepackage{graphicx}

\usepackage[T2A,T1]{fontenc}
\DeclareSymbolFont{cyrillic}{T2A}{cmr}{m}{it}
\def\makecyrsymbol#1#2{%
    \begingroup\edef\temp{\endgroup
        \noexpand\DeclareMathSymbol{\noexpand#1}
        {\noexpand\mathalpha}{cyrillic}%
        {\expandafter\expandafter\expandafter
            \calccyr\expandafter\meaning\csname T2A\string#2\endcsname\end}}%
    \temp}
\expandafter\def\expandafter\calccyr\string\char#1\end{#1}

\makecyrsymbol\Zhe\CYRZH
\makecyrsymbol\Be\CYRB
\makecyrsymbol\Shcha\CYRSHCH
\makecyrsymbol\Sha\CYRSH

\makeatletter
\def\upintkern@{\mkern-7mu\mathchoice{\mkern-3.5mu}{}{}{}}
\def\upintdots@{\mathchoice{\mkern-4mu\@cdots\mkern-4mu}%
 {{\cdotp}\mkern1.5mu{\cdotp}\mkern1.5mu{\cdotp}}%
 {{\cdotp}\mkern1mu{\cdotp}\mkern1mu{\cdotp}}%
 {{\cdotp}\mkern1mu{\cdotp}\mkern1mu{\cdotp}}}

\newcommand{\UpMultiIntegral}[1]{%
  \edef\ints@c{\noexpand\upintop
    \ifnum#1=\z@\noexpand\upintdots@\else\noexpand\upintkern@\fi
    \ifnum#1>\tw@\noexpand\upintop\noexpand\upintkern@\fi
    \ifnum#1>\thr@@\noexpand\upintop\noexpand\upintkern@\fi
    \noexpand\upintop
    \noexpand\ilimits@
  }%
  \futurelet\@let@token\ints@a
}
\makeatother

\DeclareFontFamily{OMX}{mdbch}{}
\DeclareFontShape{OMX}{mdbch}{m}{n}{ <->s * [0.8]  mdbchr7v }{}
\DeclareFontShape{OMX}{mdbch}{b}{n}{ <->s * [0.8]  mdbchb7v }{}
\DeclareFontShape{OMX}{mdbch}{bx}{n}{<->ssub * mdbch/b/n}{}

\DeclareSymbolFont{uplargesymbols}{OMX}{mdbch}{m}{n}
\SetSymbolFont{uplargesymbols}{bold}{OMX}{mdbch}{b}{n}
\DeclareMathSymbol{\upintop}{\mathop}{uplargesymbols}{82}
\DeclareMathSymbol{\upointop}{\mathop}{uplargesymbols}{"48}

\DeclareFontEncoding{MDB}{}{}
\DeclareFontFamily{MDB}{mdbch}{}
\DeclareFontShape{MDB}{mdbch}{m}{n}{ <->s * [0.8]  mdbchrmb }{}
\DeclareFontShape{MDB}{mdbch}{b}{n}{ <->s * [0.8]  mdbchbmb }{}
\DeclareFontShape{MDB}{mdbch}{bx}{n}{<->ssub * mdbch/b/n}{}
\DeclareFontSubstitution{MDB}{cmr}{m}{n}
\DeclareSymbolFont{mathdesignB}{MDB}{mdbch}{m}{n}%
\SetSymbolFont{mathdesignB}{bold}{MDB}{mdbch}{b}{n}%
\DeclareMathSymbol{\upintclockwise}{\mathop}{mathdesignB}{128}
\DeclareMathSymbol{\upointclockwise}{\mathop}{mathdesignB}{130}
\DeclareMathSymbol{\upointctrclockwise}{\mathop}{mathdesignB}{132}
\DeclareMathSymbol{\upoiint}{\mathop}{mathdesignB}{134}
\DeclareMathSymbol{\upoiiint}{\mathop}{mathdesignB}{136}

\makeatletter
\newcommand{\upint}{\DOTSI\upintop\ilimits@}
\newcommand{\upoint}{\DOTSI\upointop\ilimits@}
\makeatother

\theoremstyle{remark}

\newtheorem*{rem*}{Remark}

\renewcommand{\f}{{\varphi}}

\newcommand\dg{\frak D}

\newcommand{\ri}{{\rm i}}

\newcommand{\Bs}{\Be_{\be,1}^1}

\begin{document}

%
%



\numberwithin{equation}{section}

\numberwithin{equation}{section}

\title{Functions of pairs of unbounded noncommuting self-adjoint operators under perturbation}
\author{A.B.  Aleksandrov, V.V. Peller}
\thanks{}.


\

\begin{abstract}
For a pair $(A,B)$ of not necessarily bounded and not necessarily commuting self-adjoint operators and for a function $f$ on the Euclidean space $\R^2$ that belongs to the inhomogeneous Besov class $\Bs(\R^2)$, we define the function $f(A,B)$ of these operators as a densely defined operator.
We consider the problem of estimating the functions $f(A,B)$ under perturbations of the pair 
$(A,B)$. It is established that if 
$1\le p\le2$, and $(A_1,B_1)$ and $(A_2,B_2)$ are pairs of not necessarily bounded and not necessarily commuting self-adjoint operators such that the operators
$A_1-A_2$ and $B_1-B_2$ belong to the Schatten--von Neumann class $\bS_p$ with $p\in[1,2]$
and $f\in\Bs(\R^2)$, then the following Lipschitz type estimate holds:
$$
\|f(A_1,B_1)-f(A_2,B_2)\|_{\bS_p}
\le\const\|f\|_{\Bs}\max\big\{\|A_1-A_2\|_{\bS_p},\|B_1-B_2\|_{\bS_p}\big\}.
$$
\end{abstract} 

\maketitle


\


\

\setcounter{section}{0}
\section{\bf Introduction}
\setcounter{equation}{0}
\label{In}

\

The results of this note extend the results of \cite{ANP}
to the case of pairs of unbounded noncommuting self-adjoint operators.
Recall (see, e.g.,, \cite{ANP}), that for a pair $(A,B)$ of not necessarily bounded self-adjoint operators and for a complex-valued function $f$ on $\R\times\R$, being a Schur multiplier with respect to arbitrary Borel spectral measures, the function $f(A,B)$ of $A$ and $B$ is defined as the double operator integral 
\bay
\label{fAB}
f(A,B)\df\iint_{\R\times\R}f(x,y)\,dE_A(x)\,dE_B(y)
\df\iint_{\R\times\R}f(x,y)\,dE_A(x)I\,dE_B(y),
\ey
where $I$ is the identity operator while $E_A$ and $E_B$ are spectral measures of the operators  $A$ and $B$. Then $f(A,B)$ is a bounded operator.
We refer the reader to 
\cite{BS1}, \cite{BS2} and \cite{BS3} 
for the definition and basic properties of double operator integrals.

We also refer the reader to 
\cite{Pe1} and \cite{AP4} for the definition of Schur multipliers with respect to spectral measures. Recall (see \cite{Pe1} and \cite{AP4}) that a function $\Phi$ is a 
{\it Schur multiplier with respect to spectral measures $E_1$ and $E_2$} if and only if  $\Phi$ belongs to the  {\it Haagerup tensor product} 
$L^\be(E_1)\!\otimes_{\rm h}\!L^\be(E_2)$, i.e., $\Phi$ admits a representation of the form 
\bay
\label{Phifnpsin}
\Phi(x,y)=\sum_n\f_n(x)\psi_n(y),
\ey
where $\f_n\in L^\be(E_1)$, $\psi_n\in L^\be(E_2)$ and
\bay 
\label{koren'izproizvedeniya}
\left(\left\|\sum_n|\f_n(x)|^2\right\|_{L^\be(E_1)}
\left\|\sum_n|\psi_n(x)|^2\right\|_{L^\be(E_2)}\right)^{1/2}<\be.
\ey
The norm of $\Phi$ in $L^\be(E_1)\!\otimes_{\rm h}\!L^\be(E_2)$ is the infimum of the left-hand side of \rf{koren'izproizvedeniya} over all representations of the form
\rf{Phifnpsin}. In this case 
$$
\iint_{\X\times\Y}\Phi(x,y)\,dE_1(x)Q\,dE_2(y)=
\sum_n\left(\int\f_n\,dE_1\right)Q\left(\int\psi_n\,dE_2\right);
$$
the series on the right converges in the weak operator topology and
$$
\left\|\iint\Phi\,dE_1Q\,dE_2\right\|\le\|\Phi\|_{L^\be\!\otimes_{\rm h}\!L^\be}\|Q\|
$$
(see, e.g., \cite{AP4}).

In this note we define functions $f(A,B)$ of unbounded noncommuting operators for certain functions $f$ that do not belong to the Haagerup tensor product of the spaces of bounded functions. In this case $f(A,B)$ turns out to be a densely defined unbounded operator.

In \cite{ANP} for pairs
$(A_1,B_1)$ and $(A_2,B_2)$  of noncommuting bounded self-adjoint operators $A$ and $B$ and for functions $f$  in the {\it homogeneous Besov class} 
$B_{\be,1}^1(\R^2)$ the operators $f(A_1,B_1)$ and $f(A_2,B_2)$ were defined and
the following Lipschitz type estimate in the Schatten--von Neumann classes $\bS_p$ with $1\le p\le2$ was given:
$$
\|f(A_1,B_1)-f(A_2,B_2)\|_{\bS_p}\le\const\|f\|_{B_{\be,1}^1}
\max\big\{\|A_1-B_1\|_{\bS_p},\|A_2-B_2\|_{\bS_p}\big\}.
$$
In the same paper \cite{ANP} it was shown that the same inequality is false for
$p>2$ and is false in the operator norm.

Recall also that in the case of functions of a single self-adjoint operator, such Lipschitz type estimates hold true for $1\le p\le\be$, see \cite{Pe1} and \cite{Pe2}. 

The main purpose of this note is to establish this inequality for pairs of unbounded noncommuting self-adjoint operators and for functions$ f$ in the {\it inhomogeneous} Besov class $\Bs(\R^2)$. We refer the reader to \cite{Pee} for the definition and basic properties of Besov spaces.

As in the case of bounded noncommuting operators, the key role is played by triple operator integrals. We refer the reader to
\cite{ANP} and \cite{AP5}.

\

\section{\bf Triple operator integrals, Haagerup and Haagerup-like tensor products}
\setcounter{equation}{0}
\label{TrOi}

\

Triple operator integrals are expressions of the form
\bay
\label{troinoi}
\int\limits_{\X_1}\int\limits_{\X_2}\int\limits_{\X_3} 
\Psi(x_1,x_2,x_3)\,dE_1(x_1)T\,dE_2(x_2)R\,dE_3(x_3),
\ey
where $\Psi$ is a bounded measurable function on $\X_1\times\X_2\times\X_3$; 
$E_1$, $E_2$ are $E_3$ spectral measures on Hilbert space while $T$ and $R$ are bounded  linear operators.
Such operator integrals can be defined under certain assumptions on $\Psi$, $T$ and $R$.

In \cite{Pe3} integrals of the form \rf{troinoi} are defined for arbitrary bounded operators $T$ and $R$ and for functions $\Psi$ in the integral projective tensor product 
$L^\be(E_1)\!\otimes_{\rm i}\!L^\be(E_2)\!\otimes_{\rm i}\!L^\be(E_3)$. In this case the following inequality holds:
 $$
 \left\|\iiint\Psi\,dE_1T\,dE_2)R\,dE_3\right\|_{\bS_r}\le
 \|\Psi\|_{L^\be\otimes_{\rm i}L^\be\otimes_{\rm i}L^\be}\|T\|_{\bS_p}\|R\|_{\bS_q},\quad
 T\in\bS_p,~R\in\bS_q,
 $$
 $$
 \frac1r=\frac1p+\frac1q\quad\mbox{under the assumption}\quad\frac1p+\frac1q\le1.
 $$

Later in \cite{JTT} triple operator integrals were defined for functions 
$\Psi$ in the {\it Haagerup tensor product} 
$L^\be(E_1)\!\otimes_{\rm h}\!L^\be(E_2)\!\otimes_{\rm h}\!L^\be(E_3)$. We refer the reader to
\cite{AP5} for the definition and basic properties of such Haagerup tensor products. Note here that for functions $\Psi$ in 
$L^\be(E_1)\!\otimes_{\rm h}\!L^\be(E_2)\!\otimes_{\rm h}\!L^\be(E_3)$, the following estimates hold: 
$$
\left\|\iiint\Psi\,dE_1T\,dE_2R\,dE_3\right\|
\le\|\Psi\|_{L^\be\!\otimes_{\rm h}\!L^\be\!\otimes_{\rm h}\!L^\be}
\|T\|\cdot\|R\|
$$
in the case of bounded operators $T$ and $R$, and
$$
\left\|\iiint\!\Psi\,dE_1T\,dE_2R\,dE_3\right\|_{\bS_r}
\!\!\!\le\|\Psi\|_{L^\be\!\otimes_{\rm h}\!L^\be\!\otimes_{\rm h}\!L^\be}
\|T\|_{\bS_p}\|R\|_{\bS_q}
$$
in the case when $T\in\bS_p$, $R\in\bS_q$, $1/r=1/p+1/q$ and $p,\,q\in[2,\be]$.

However, it turned out that for Lipschitz type estimates for functions of pairs of noncommuting operators, we need triple operator integrals with integrands in so-called Haagerup-like tensor products of the first kind $L^\be(E_1)\!\otimes_{\rm h}\!L^\be(E_2)\!\otimes^{\rm h}\!L^\be(E_3)$
and the second kind $L^\be(E_1)\!\otimes^{\rm h}\!L^\be(E_2)\!\otimes_{\rm h}\!L^\be(E_3)$.
Such tensor products were introduced in \cite{ANP} and were studied in more detail in \cite{AP5}.

In \cite{ANP} and \cite{AP5} it was shown that if
$\Psi\in L^\be(E_1)\!\otimes_{\rm h}\!L^\be(E_2)\!\otimes^{\rm h}\!L^\be(E_3)$, $1\le p\le2$, 
$T\in\bS_p$d an $R$ is a bounded linear operator, one can define the triple operator integral of the form \rf{troinoi}; moreover, the following estimate holds:
$$
\left\|
\iiint 
\Psi\,dE_1T\,dE_2R\,dE_3
\right\|_{\bS_p}
\le\|\Psi\|_{L^\be\otimes_{\rm h}\!L^\be\otimes^{\rm h}\!L^\be}
\|T\|_{\bS_p}\|R\|,\quad1\le p\le2.
$$
On the other hand, if
$\Psi\in L^\be(E_1)\!\otimes^{\rm h}\!L^\be(E_2)\!\otimes_{\rm h}\!L^\be(E_3)$,
$1\le p\le2$, $T$ is a bounded linear operator and $R\in\bS_p$, then
$$
\left\|
\iiint 
\Psi\,dE_1T\,dE_2R\,dE_3
\right\|_{\bS_p}\le
\|\Psi\|_{L^\be\otimes^{\rm h}\!L^\be\otimes_{\rm h}\!L^\be}
\|T\|\cdot\|R\|_{\bS_p}.
$$

Moreover,
\bay
\label{eshchorazW}
\|W\|_{\bS_p}\le
\|\Psi\|_{L^\be\otimes^{\rm h}\!L^\be\otimes_{\rm h}\!L^\be}
\|T\|\cdot\|R\|_{\bS_p}.
\ey

Note that in \cite{ANP} more general Schatten--von Neumann estimates for triple operator integrals
were obtained in the case when the integrand belongs to Haagerup-like tensor products of $L^\be$ spaces.
Later in \cite{AP5} the estimates obtained in \cite{ANP} were extended to even a more general case. 

Note also that in the same way one can define Haagerup-like tensor products  
$\B^\be\!\otimes_{\rm h}\!\B^\be\!\otimes^{\rm h}\!\B^\be$ and 
$\B^\be\!\otimes^{\rm h}\!\B^\be\!\otimes_{\rm h}\!\B^\be$, where $\B^\be$ is the space of bounded Borel functions on $\R$.

Consider now a continuously differentiable function $f$ on $\R^2$ and define the divided differences
$\dg^{[1]}f$ and $\dg^{[2]}f$ by
$$
\big(\dg^{[1]}f\big)(x_1,x_2,y)\df\frac{f(x_1,y)-f(x_2,y)}{x_1-x_2},
\quad x_1\ne x_2
$$
and
$$
\big(\dg^{[2]}f\big)(x,y_1,y_2)\df\frac{f(x,y_1)-f(x,y_2)}{y_1-y_2},\quad y_1\ne y_2.
$$
In the case when $x_1=x_2$ or $y_1=y_2$, in the definition of $\dg^{[1]}f$ and 
$\dg^{[2]}f$ one has to replace divided differences with the corresponding partial derivatives. 

We define the class $\E_\s^\be(\R^2)$ for $\s>0$ as follows:
$$
\E_\s^\be(\R^2)\df\Big\{f\in L^\be(\R^2):~\supp\F f\subset\{t\in\R^2:~\|t\|_2\le\s\}\Big\};
$$
here we use the notation $\F$ for Fourier transform.

It was established in \cite{ANP} that for $\s>0$ and $f\in\E_\s^\be(\R^2)$ the following estimates hold: 
\bay
\label{otsenkaD1}
\big\|\dg^{[1]}f\big\|_{\B^\be\otimes_{\rm h}\B^\be\otimes^{\rm h}\B^\be}
\le\const\s\|f\|_{L^\be(\R^2)},
\ey
\bay
\label{otsenkaD2}
\big\|\dg^{[2]}f\big\|_{\B^\be\otimes^{\rm h}\B^\be\otimes_{\rm h}\B^\be}
\le\const\s\|f\|_{L^\be(\R^2)}.
\ey
It follows that if a function $f$ belongs to the homogeneous Besov class
$B_{\be,1}^1(\R^2)$, then $\dg^{[1]}f\in\B^\be\!\otimes_{\rm h}\!\B^\be\!\otimes^{\rm h}\!\B^\be$ and
$\dg^{[2]}f\in\B^\be\!\otimes^{\rm h}\!\B^\be\!\otimes_{\rm h}\!\B^\be$; moreover,
$$
\big\|\dg^{[1]}f\big\|_{\B^\be\otimes_{\rm h}\B^\be\otimes^{\rm h}\B^\be}
\le\const\|f\|_{B_{\be,1}^1},
$$
$$
\big\|\dg^{[2]}f\big\|_{\B^\be\otimes^{\rm h}\B^\be\otimes_{\rm h}\B^\be}
\le\const\|f\|_{B_{\be,1}^1}.
$$

The main result of \cite{ANP} is the following fact:

\medskip

{\bf Theorem 2.1.} {\it Let $1\le p\le2$, $f\in B_{\be,1}^1(\R^2)$, and $(A_1,B_1)$ and $(A_2,B_2)$ are pairs of bounded noncommuting self-adjoint operators such that $A_2-A_1\in\bS_p$ and $B_2-B_1\in\bS_p$. Then
\begin{align*}
f(A_1,B_1)-f(A_2,B_2)
&=\!
\iiint\frac{f(x_1,y)-f(x_2,y)}{x_1-x_2}
dE_{A_1}(x_1)(A_1-A_2)dE_{A_2}(x_2)dE_{B_1}(y)\nonumber\\[.2cm]
&+\iiint\frac{f(x,y_1)-f(x,y_2)}{y_1-y_2}
dE_{A_2}(x)dE_{B_1}(y_1)(B_1-B_2)dE_{B_2}(y_2).
\end{align*}
Moreover, the following estimate holds:
$$
\|f(A_1,B_1)-f(A_2,B_2)\|_{\bS_p}\le\const\|f\|_{B_{\be,1}^1}
\max\big\{\|A_1-A_2\|_{\bS_p},\|B_1-B_2\|_{\bS_p}\big\}.
$$
}

The main objective of this note is to establish the same inequality in the case of unbounded noncommuting pairs of operators under the assumption that the function $f$ 
belongs to the inhomogeneous Besov class  $\Bs(\R^2)$.

\

\section{\bf Functions of pairs of unbounded noncommuting self-adjoint operators}
\label{nekommu}

\

Recall that we have defined functions of not necessarily commuting self-adjoint operators by  \rf{fAB} in the case when thefunction $f$ belongs to the Haagerup tensor product 
$\B^\be\!\otimes_{\rm h}\!\B^\be$. Moreover, the following estimate holds: 
$$
\|f(A,B)\|\le\|f\|_{\B^\be\otimes_{\rm h}\B^\be},\quad f\in\B^\be\otimes_{\rm h}\B^\be.
$$

Let $f$ be a function of two variables and let $f_\sharp$ be the function defined by the equality
$f_\sharp(s,t)\df(1-\ri t)^{-1}f(s,t)$.
Suppose that $f_\sharp\in\B^\be\!\otimes_{\rm h}\!\B^\be$. 
We define the operator $f(A,B)$ by
$$
f(A,B)\df f_\sharp(A,B)(I-\ri B)=\left(\,\,\,\iint\limits_{\R\times\R}f_\sharp(s,t)\,dE_A(s)\,dE_B(t)\right)(I-\ri B).
$$
Then $f(A,B)$ is a densely defined operator whose domain coincides with the domain $D(B)$ of $B$. It does not have to be bounded but the operator $f(A,B)(I-\ri B)^{-1}$ is bounded.

Note that if $f\in\E_\s^\be(\R^2)$, $\s>0$, then
 $f_\sharp\in\B^\be\!\otimes_{\rm h}\!\B^\be$. This was established in Corollary 7.3 of \cite{AP6} for functions $f$ in $\big(\E_\s^\be(\R^2)\big)_+$, i.e., functions $f$ in $\E_\s^\be(\R^2)$ whose Fourier transform is supported in $[0,\be)\times[0,\be)$. Clearly, the same is also true for  $f\in\E_\s^\be(\R^2)$, $\s>0$.
Thus, if $f\in\E_\s^\be(\R^2)$, then
the operator $f_\sharp(A,B)$ is bounded, while 
$f(A,B)$ is a not necessarily bounded densely defined operator with domain $D(B)$. Moreover,
$$
\|f_\sharp\|_{\B^\be\otimes_{\rm h}\B^\be}\le\const(1+\s)\|f\|_{L^\be(\R^2)}.
$$
This can be verified in the same was as in Corollary 7.3 of \cite{AP6}.

\begin{thm} 
\label{th73}
Let $f\in\Bs(\R^2)$. Then $f_\sharp\in \B^\be\otimes_{\rm h}\B^\be$
and $\|f_\sharp\|_{\B^\be\otimes_{\rm h}\B^\be}\le\const\|f\|_{\Bs}$.
\end{thm}

\

\section{\bf Integral formulae for operator differences and Lipschitz type estimates}
\label{OtsLiptipa}

\

In this section we state the main result of the note. We obtain a formula for the operator difference in terms of triple operator integrals and we establish a Lipschitz type estimate in the 
$\bS_p$ norm for $p\in[1,2]$. We are going to deal with not necessarily bounded and not necessarily commuting self-adjoint operators.

\begin{thm}
\label{teor1}
Let $f\in \E^\be_\s(\R^2)$, and $A_1$, $A_2$ and $B$ are self-adjoint operators such that
$A_1-A_2\in\bS_2$. Then
$$
f(A_1,B)-f(A_2,B)=\iiint\frac{f(x_1,y)-f(x_2,y)}{x_1-x_2}
\,dE_{A_1}(x_1)(A_1-A_2)\,dE_{A_2}(x_2)\,dE_{B}(y),
$$
and so
$$
\|f(A_1,B)-f(A_2,B)\|_{\bS_p}\le\const\s\|f\|_{L^\be(\R^2)}\|A_1-A_2\|_{\bS_p}.
$$
\end{thm}

Recall that $\dg^{[1]}f\in\B^\be\otimes_{\rm h}\B^\be\otimes^{\rm h}\B^\be$ (see \rf{otsenkaD1}), and so the triple operator integral on the right is defined.

\begin{cor} 
\label{cor58}
Let $f\in \Bs(\R^2)$ and $1\le p\le2$. Suppose that $A_1$, $A_2$ and $B$ are self-adjoint operators such that 
$A_2-A_1\in\bS_p$. Then the following inequality holds:
$$
\|f(A_1,B)-f(A_2,B)\|_{\bS_p}\le\const\|f\|_{\Bs}\|A_1-A_2\|_{\bS_p}.
$$
\end{cor}

\begin{thm}
\label{teor2}
Let $f\in \E^\be_\s(\R^2)$. Suppose that $A$, $B_1$ and $B_2$ are
self-adjoint operators such that 
$B_2-B_1\in\bS_2$. Then the following equality holds:
$$
f(A,B_1)-f(A,B_2)=
\iiint\frac{f(x,y_1)-f(x,y_2)}{y_1-y_2}
\,dE_{A}(x)\,dE_{B_1}(y_1)(B_1-B_2)\,dE_{B_2}(y_2).
$$
\end{thm}

Again, $\dg^{[2]}f\in\B^\be\otimes^{\rm h}\B^\be\otimes_{\rm h}\B^\be$ (see \rf{otsenkaD2}), and so the triple operator integral on the right is defined.

\begin{cor}
\label{cor510}
Let $f\in \E^\be_\s(\R^2)$ with $p\in[1,2]$. Suppose that $A$, $B_1$ and $B_2$ are
self-adjoint operators such that
$B_2-B_1\in\bS_p$. Then
$$
\|f(A,B_1)-f(A,B_2)\|_{\bS_p}\le\const\s\|f\|_{L^\be(\R^2)}\|B_1-B_2\|_{\bS_p}.
$$
\end{cor}

\begin{thm}
\label{lipschitseva_otsenka_dlya_Besova}
Let $f\in\Bs(\R^2)$ with $p\in[1,2]$. Suppose that $A_1$, $A_2$, $B_1$ and $B_2$ are
self-adjoint operators such that
$A_2-A_1\in\bS_p$ and $B_2-B_1\in\bS_p$. Then
$$
\|f(A_1,B_1)-f(A_2,B_2)\|_{\bS_p}\le\const\|f\|_{\Bs}
\max\big\{\|A_1-A_2\|_{\bS_p},\|B_1-B_2\|_{\bS_p}\big\}.
$$
\end{thm}

\begin{thm} 
\label{ab12}
Let $f\in\Bs(\R^2)$. Suppose that $A_1$, $A_2$, $B_1$ and $B_2$ are
self-adjoint operators such that 
$A_1-A_2\in\bS_2$ and $B_1-B_2\in\bS_2$. Then the following identity holds:
\begin{align*}
f(A_1,B_1)&-f(A_2,B_2)\nonumber\\[.2cm]
&=
\iiint\frac{f(x_1,y)-f(x_2,y)}{x_1-x_2}
\,dE_{A_1}(x_1)(A_1-A_2)\,dE_{A_2}(x_2)\,dE_{B_1}(y)\\[.2cm]
&+\iiint\frac{f(x,y_1)-f(x,y_2)}{y_1-y_2}
\,dE_{A_2}(x)\,dE_{B_1}(y_1)(B_1-B_2)\,dE_{B_2}(y_2).
\end{align*}
\end{thm}

%
%

\

\

\section*{SOURCE OF FUNDING}

\

The research on \S\,2  is supported by Russian Science Foundation [grant number
18-11-00053]. 
The research on \S\,3 is supported by Russian Science Foundation [grant number 20-61-46016].
The remaining results are supported by
a grant of the Government of the Russian Federation for the state support of scientific research, carried out under the supervision of leading scientists, agreement  075-15-2021-602

\
 
 \begin{footnotesize}
 
\noindent
\begin{tabular}{p{7cm}p{15cm}}
A.B. Aleksandrov & V.V. Peller \\
St.Petersburg Department & Department of Mathematics\\
Steklov Institute of Mathematics  & and Computer Sciences\\
Fontanka 27, 191023 St.Petersburg & St.Petersburg State University\\
Russia&Universitetskaya nab., 7/9,\\
email: alex@pdmi.ras.ru&199034 St.Petersburg, Russia\\
\\
&St.Petersburg Department\\
&Steklov Institute of Mathematics\\
&Russian Academy of Sciences\\
&Fontanka 27, 191023 St.Petersburg\\
&Russia\\
\end{tabular}
\end{footnotesize}

\


\begin{thebibliography}{99}
%

%

\bibitem{ANP} {\sc A.B. Aleksandrov, F.L. Nazarov} and
{\sc V.V. Peller}, {\em Functions of noncommuting self-adjoint operators under perturbation and estimates of triple operator integrals},  Adv. Math. {\bf295} (2016), 1--52.



%


\bibitem{AP4} {\sc A.B. Aleksandrov} and {\sc V.V. Peller}, {\it Operator Lipschitz functions}, 
Uspekhi Matem. Nauk {\bf71:4} (2016), 3--106.

English transl.: Russian Mathematical Surveys {\bf71:4} (2016), 605--702.

\bibitem{AP5} {\sc A.B. Aleksandrov} and {\sc V.V. Peller}, {\it Multiple operator integrals, Haagerup and Haagerup-like tensor products, and operator ideals}, Bulletin London Math. Soc. {\bf49} (2016), 463--479.



\bibitem{AP6}  {\sc A.B. Aleksandrov} and {\sc V.V. Peller},  {\em Functions of perturbed commuting dissipative operators},
Math. Nachr., to appear

%

\bibitem{BS1} {\sc M.S. Birman} and {\sc M.Z. Solomyak},
{\em Double Stieltjes operator integrals},
Problems of Math. Phys., Leningrad. Univ. {\bf1} (1966), 33--67 (Russian).
\newline
English transl.: Topics Math. Physics {\bf1} (1967), 25--54, Consultants Bureau Plenum
Publishing Corporation, New York.

\bibitem{BS2} {\sc M.S. Birman} and {\sc M.Z. Solomyak},
 {\em Double Stieltjes operator integrals. II},
 Problems of Math. Phys., Leningrad. Univ. {\bf2} (1967), 26--60 (Russian).
 \newline
English transl.: Topics Math. Physics {\bf2} (1968), 19--46, Consultants Bureau Plenum
Publishing Corporation, New York.

\bibitem{BS3} {\sc M.S. Birman} and {\sc M.Z. Solomyak},
{\em Double Stieltjes operator integrals. III},
Problems of Math. Phys., Leningrad. Univ. {\bf6} (1973), 27--53 (Russian).

%







%



%
%

\bibitem{JTT} {\sc K. Juschenko, I.G. Todorov} and {\sc L. Turowska}, {\em Multidimensional operator multipliers}, Trans. Amer. Math. Soc. {\bf361}
(2009), 4683--4720.

\bibitem{Pee} {\sc J. Peetre},
{\em New thoughts on Besov spaces}, Duke Univ. Press., Durham, NC, 1976.

\bibitem{Pe1} {\sc V.V. Peller},
{\em Hankel operators in the theory of perturbations of unitary and self-adjoint operators},
Funktsional. Anal. i Prilozhen. {\bf19:2}  (1985),
37--51 (Russian).

English transl.: Funct. Anal. Appl. {\bf19} (1985), 111--123.


\bibitem{Pe2} {\sc V.V. Peller},
{\em Hankel operators in the perturbation theory of unbounded self-adjoint operators},
Analysis and partial differential equations,  529--544,
Lecture Notes in Pure and Appl. Math., {\bf122}, Marcel Dekker, New York, 1990.


\bibitem{Pe3} {\sc V.V. Peller}, {\it Multiple operator integrals in perturbation theory},  Bull. Math. Sci. {\bf6} (2016), 15--88.




%

 
\end{thebibliography}
\end{document}